\input amstex
\input epsf
\documentstyle{amsppt}
\magnification 1200
\vcorrection{-9mm}
\NoBlackBoxes

\def\C{\Bbb C}

\def\CP{\Bbb{CP}}
\def\eps{\varepsilon}
\def\Arg{\operatorname{Arg}}
\def\pr{\operatorname{pr}}
\def\op{\operatorname{op}}

\def\ord{\operatorname{ord}}

\def\QP{\Cal Q}

\def\SCB{\Cal{SB}}

\def\CB{\Cal B}

\def\Im{\operatorname{Im}}
\def\lk{\operatorname{lk}}

\def\refBO     {1}
\def\refBR     {2}
\def\refEN     {3}
\def\refFS     {4}
\def\refFW     {5}
\def\refGO     {6}
\def\refHay    {7}
\def\refKM     {8}
\def\refLMknot {9}
\def\refLMlink {10}
\def\refM      {11}
\def\refMo     {12}
\def\refN      {13}
\def\refNW     {14}
\def\refO      {15}
\def\refRuTop  {16}
\def\refRuFancy{17}

\def\figSym    {1}
\def\figHopf   {2} \let\figLink=\figHopf
\def\figENa    {3}
\def\figENb    {4}
\def\figFiveCrA{5}
\def\figTorus  {6}
\def\figFiveCr {7}

\def\eqA   {1}

\def\sectSimpleEx {2}
\def\sectBR       {3} \let\sectKM=\sectBR
\def\sectHopf     {4}
\def\sectFurther  {5}
\def\sectFiveCr   {6}

\def\thH   {1.1}
\def\thO   {1.2}

\def\thA   {\sectSimpleEx.1}
\def\remA  {\sectSimpleEx.2}
\def\corA  {\sectSimpleEx.3}
\def\thB   {\sectSimpleEx.4}

\def\thKM   {\sectKM.1}
\def\propBR {\sectKM.2}
\def\remBR  {\sectKM.3}
\def\remM   {\sectKM.4}
\def\defOut {\sectKM.5}

\def\propSum{\sectKM.6}
\def\corBR  {\sectKM.7}
\def\lemBBS {\sectKM.8}

\def\remWermer  {\sectHopf.1}

\def\defIterTor {\sectFurther.1}
\def\remIter    {\sectFurther.2}
\def\propIter   {\sectFurther.3}
\def\lemIter    {\sectFurther.4}

\def\thFW       {\sectFiveCr.1}
\def\propCaseH  {\sectFiveCr.2}
\def\propFiveCr {\sectFiveCr.3}

\topmatter
\title  Plane algebraic curves in fancy balls
\endtitle

\author
        N.~G.~Kruzhilin and S.~Yu.~Orevkov
\endauthor

\address
        Steklov Mathematical Institute of Russian Academy of Sciences (N.K. and S.O.)
\endaddress

\address
        IMT, l'universit\'e Paul Sabatier, 118 route de Narbonne, Toulouse, France (S.O.)
\endaddress


\abstract
        Boileau and Rudolph [\refBR] called a link $L$ in the $3$-sphere
        a {\it $\C$-boundary} if it can be realized as the intersection
        of an algebraic curve $A$ in $\C^2$ with the boundary of a
        smooth embedded $4$-ball $B$. They showed that some links are
        not $\C$-boundaries. We say that $L$ is a {\it strong
        $\C$-boundary} if $A\setminus B$ is connected. In particular, all
        quasipositive links are strong $\C$-boundaries.

        In this paper we give examples of non-quasipositive strong $\C$-boundaries
        and non-strong $\C$-boundaries. We give a complete classification
        of (strong) $\C$-boundaries with at most 5 crossings.
\endabstract

\endtopmatter

\document

\head   \S1. Introduction   \endhead
Let $B\subset\C^2$ be diffeomorphic to a $4$-ball
and $A$ be a complex analytic curve in a neighbourhood of $B$ which is transverse to
$\partial B$ (since we consider only topological properties, we may assume that
$A$ is a piece of an algebraic curve).
Let $L=A\cap\partial B$ be a link in the $3$-sphere $\partial B$,
endowed with the boundary orientation from $A\cap B$.
Which links can be obtained in this way?
(All links in this paper are assumed to be oriented.)

If we impose no additional restrictions then, as Lee Rudolph  showed
in [\refRuFancy],%
\footnote{The title of our paper is a rephrasing of the title of [\refRuFancy].}
the answer is ``any link''. Moreover, any oriented surface without
closed components can be realized as $A\cap B$.

If $B$ is strictly pseudoconvex (for example, a usual round ball), then
it was shown in [\refBO] that a link is realizable in this way if and only if
it is {\it quasipositive} (the ``if'' part being earlier proven in [\refRuTop]),
i.e., it is the braid closure of a quasipositive braid
(an $n$-braid is called {\it quasipositive} if it is a product of conjugates
of the standard generators $\sigma_1,\dots,\sigma_{n-1}$ of the
braid group $B_n$). This is a rather strong restriction on the class of
possible links (see [\refBR, \refHay]).

In [\refBR], a link is called a {\it $\C$-boundary} if it is
realizable as $A\cap\partial B$ where $B$ is diffeomorphic to a 4-ball (without
any pseudoconvexity assumptions) and $A$ is a whole algebraic curve in $\C^2$,
not just a piece of an algebraic curve as in [\refRuFancy].
It is also natural to distinguish the case when $L$ is realizable as $A\cap\partial B$
as above, and moreover $A\setminus B$ is connected. We call such links
{\it strong $\C$-boundaries}.
It was observed in [\refBR] that Kronheimer and Mrowka's result [\refKM]
implies some restrictions on this class of links, in particular, there exist
knots and links that are not concordant to any $\C$-boundary.
In fact, some results stated in [\refBR] for arbitrary $\C$-boundaries are true only for
strong $\C$-boundaries; see more details in \S\sectBR.

Michel Boileau (private communication) asked if there exist non-quasipositive
$\C$-boundaries.
Here we give an affirmative answer to this question. Moreover, we show
that all the following inclusions are strict:
$$
\split
   \QP:=\{\text{quasipositive links}\}\subset
   \SCB&:=\{\text{strong $\C$-boundaries}\}\\
   &\subset
   \CB:=\{\text{$\C$-boundaries}\}\subset
   %\LL:=
   \{\text{all links}\}.
\endsplit
$$
To prove that some $\C$-boundaries are not quasipositive,
we use the following facts.

\proclaim{ Theorem \thH } {\rm([\refHay, Cor.~1.5]).}
If a link and its mirror are both quasipositive, then the link is trivial.
\endproclaim

\proclaim{ Theorem \thO } {\rm([\refO, Thms.~1.1 and 1.2]).}
If the split sum or a connected sum of links $L_1$ and $L_2$ is quasipositive,
then $L_1$ and $L_2$ are quasipositive links.
\endproclaim

Another necessary condition for the quasipositivity of links follows from
the Franks--Williams--Morton inequality (see Theorem~\thFW\ below).

The definition of strong $\C$-boundary can be reformulated equivalently by
replacing the condition that $A\setminus B$ is connected by the condition
that $A\setminus B$ does not have bounded components.
If $B$ is strictly pseudoconvex,  bounded components may appear:
see Wermer's example [\refNW, p.~34] (but no component of $A\setminus B$ cannot
be a disk by Nemirovski's result in [\refN]). Nonetheless, when $B$ is strictly
pseudoconvex, $A\cap B$ is a strong $\C$-boundary, because it is
a quasipositive link by [\refBO], thus it can be realized on the
standard round sphere by [\refRuTop], and the absence of bounded components in that case
follows from the maximum principle. So, we have $\QP\subset\SCB$.
Notice also that  Wermer's example
also provides
a non-quasipositive $\C$-boundary; see Remark~\remWermer.
%\S\sectHopf.

\smallskip
\subhead Plan of the paper \endsubhead
In \S2\ we give simplest examples of non-quasipositive $\C$-boundaries.
They are obtained by choosing a ``fancy $4$-ball'' which is a small thickening
of a $3$-ball embedded in the standard $3$-sphere.

In \S3\ we present some tools to prove that some links are not (strong) $\C$-boundaries.
All of them are based on the Kronheimer--Mrowka Theorem.

In \S\S4--5\ we discuss some links which are cut by a complex line on
an embedded $3$-sphere. If $L$ is such a link, then both $L$ and its
mirror $L^*$ are $\C$-boundaries, thus one of them is a non-quasipositive
$\C$-boundary by Theorem~\thH. In \S5\ we show that these links are iterated torus links
and we describe their Eisenbud--Neumann splice diagrams.

In \S6\ we give a complete classification of $\C$-boundaries and strong $\C$-boundaries with
up to five crossings, which shows, in particular, that all the inclusions
$\QP\subset\SCB\subset\CB\subset\{$all links$\}$ are strict. This classification easily
follows from the general facts established in the previous sections, except for the $\C$-boundary
realization of the link $5_1^2$ which is a little bit tricky.

%====================================================

\head   \S\sectSimpleEx.
         The simplest examples of non-quasipositive $\C$-boundaries
\endhead

For a link $L$, let $L^*$ denote its mirror image and
let $-L$ denote $L$ with the opposite orientation.

\proclaim{ Theorem \thA }
Let $B$ and $B_0$ be $4$-balls smoothly embedded in $\C^2$, and let
$\overline B_0\subset\operatorname{Int}B$ be a
 ball contained in $B$. Let $A$ be an algebraic curve in $\C^2$
which is transverse to $\partial B$ and $\partial B_0$.
Let $L$ and $L_0$ be the links cut by $A$ on $\partial B$ and $\partial B_0$, respectively.
Then the split sum
$L\sqcup(-L_0^*)$ and a connected sum $L\#(-L_0^*)$ (see Remark~\remA)
are $\C$-boundaries.

If $B_0$ is, moreover,  strictly pseudoconvex and
$L_0$ is non-trivial, then $L\sqcup(-L_0^*)$ and $L\#(-L_0^*)$
are non-quasipositive $\C$-boundaries.
\endproclaim

\noindent{\bf Remark \remA. } A connected sum of two links $L=L_1\# L_2$
usually depends on the choice of the components
which are joined into a single component of $L$.
In Theorem~\thA, the components
of $A\cap\partial B$ and $A\cap\partial B_0$ that we choose should adjoin
   the same connected component of $A\cap(B\setminus B_0)$.
\smallskip

\demo{ Proof }
Let us show that the link under  discussion is a $\C$-boundary. Indeed,
let $I$ be an embedded line segment in $B\setminus B_0$ which connects
$A\cap\partial B$ with $A\cap\partial B_0$. Let $U$ be a small tubular
neighbourhood of $I$. Then $B\setminus(B_0\cup U)$ is a 4-ball, and
the link cut on it by $A$ is $L\sqcup(-L_0^*)$ (if $I$ is disjoint from  $A$)
or $L\#(-L_0^*)$ (if $I\subset A$).

If  $B_0$ is strictly pseudoconvex, then $L_0$ is a quasipositive link.
If it is non-trivial, Theorem~\thH\ implies that
$-L_0^*$ is not quasipositive and the result follows from Theorem~\thO.
\qed\enddemo

\proclaim{ Corollary \corA }
Let $L$ be a non-trivial quasipositive link. Then
$L\sqcup(-L^*)$ and $L\#(-L^*)$ are non-quasipositive $\C$-boundaries.
If, moreover, $L$ is a knot, then $L\#(-L^*)$ is a non-quasipositive
strong $\C$-boundary.
\endproclaim

This construction admits the following generalization.

\proclaim{ Theorem \thB }
Let $L$ be a $\C$-boundary link in $S^3$  transverse to
an equatorial $2$-sphere $S^2\subset S^3$. Let $H$ be one
of the halves of $S^3\setminus S^2$ and $\xi:S^3\to S^3$ be
the symmetry with respect to $S^2$.
Then the link $(L\cap H)\cup\xi(-L\cap H)$ is a
non-quasipositive $\C$-boundary unless it is a trivial link.
\endproclaim

\midinsert
\centerline{\epsfxsize=45mm\epsfbox{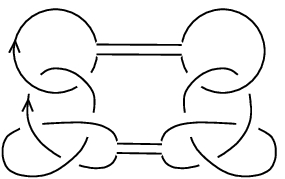}}
\botcaption{Figure~\figSym}
L10n36(1) in [\refLMlink].
\endcaption
\endinsert

\demo{ Proof }
Let $(A,B)$ be a $\C$-boundary realization of $L$. Let $f:S^3\to\partial B$
be a diffeomorphism that maps $L$ to $A\cap\partial B$.
Let $B'$ be a small thickening of $f(H)$. Then $(A,B')$ is
a $\C$-boundary realization of the required link.
Being amphicheiral, it is either trivial or non-quasipositive by Theorem~\thH.
\qed\enddemo

This theorem allows us to
construct a lot of non-quasipositive $\C$-boundaries.
In Figure~\figSym\ we give an example of such a link.
Starting with any quasipositive braid, one can construct many others.

\head
\S\sectBR. Restrictions on (strong) $\C$-boundaries
\endhead

All the known restrictions on (strong) $\C$-boundaries are more or less
immediate consequences of the Kronheimer--Mrowka theorem [\refKM] (also known
as the Thom Conjecture or the Adjunction Inequality) and its version for
immersed 2-surfaces in $\CP^2$ with negative double points
(see [\refFS], [\refM, \S2]), which was actually proven in [\refKM]
but was not explicitly formulated there.

\proclaim{ Theorem \thKM } {\rm(the Immersed Thom Conjecture).}
Let $\Sigma$ be a connected oriented closed surface of genus $g$ and
$j:\Sigma\to\CP^2$
be an immersion which has only negative ordinary double points as self-crossings.
Let $j_*([\Sigma])=d[\CP^1]\in H_2(\CP^2)$ with $d>0$. Then $g$ is bounded below by the
genus of a smooth algebraic curve of degree $d$, that is,
$g \ge (d-1)(d-2)/2$.
\endproclaim

Given a link $L$ in $S^3=\partial B^4$, we define the {\it slice
Euler characteristic} of $L$ by $\chi_s(L)=\max_\Sigma\chi(\Sigma)$ where the
maximum is taken over all embedded smooth oriented surfaces $\Sigma$ without
closed components, such that $\partial\Sigma=L$.

Similarly, we define the {\it slice negatively immersed
Euler characteristic} of $L$ by $\chi_s^-(L)=\max_{(\Sigma,j)}\chi(\Sigma)$ where the
maximum is taken over all immersions $j:(\Sigma,\partial\Sigma)\to(B^4,S^3)$
of oriented surfaces $\Sigma$ without closed components such that
$j(\Sigma)$ has only negative double points and $\partial\Sigma=L$.

Theorem \thKM\ immediately implies the following.

\proclaim{ Proposition \propBR } {\rm(cf. [\refBR, Thm.~1.3])}
Let $A$ be a smooth algebraic curve in $\C^2$ which is  transverse to the boundary of
a $4$-ball $B$  smoothly embedded in $\C^2$ and such  that
$A\setminus B$ {\bf is connected},
and let $L=A\cap\partial B$. Then $\chi_s^-(L)=\chi_s(L)=\chi(A\cap B)$.
\endproclaim

The connectedness condition in Proposition~\propBR\ can be replaced by
the condition that $A\setminus B$ does not have bounded components.
Indeed, in this case $A\setminus B$ becomes connected after a perturbation
of the union of $A$ and a generic line, far from $B$.

\demo{Proof} We replace $A\setminus B$ by a negatively immersed
surface $j(\Sigma)$ of maximum $\chi(\Sigma)$ and apply Theorem~\thKM.
\qed\enddemo

\noindent{\bf Remark \remBR.} The connectedness condition is missing
in [\refBR, Thm.~1.3]. Without this condition, Proposition~\propBR\ is wrong.
Indeed, let $A=\{w=0\}$ and let $B$ be the unit ball ``drilled''
along the line segment $[(0, 0), (0, 1)]$. Then $\chi(A\cap B)=0$ whereas
$\chi_s(L)=2$. The proof fails because, when we replace $A\cap B$ with $\Sigma$,
 the Euler characteristic increases  due to  splitting out  a
  2-sphere, while the Euler characteristic of the unbounded component does
not change.
Note that [\refBR, Prop.~1.4] is wrong even for strong $\C$-boundaries if
both links are multi-component. A correct version is Proposition~\propSum\ below.
\medskip

\noindent{\bf Remark \remM.}
A similar inaccuracy appears in [\refM]: it has not been checked if the auxiliary surface
(analog of $(B\setminus A)\cup\Sigma$ in the proof of Theorem~\thKM) is connected.
So the conclusion of [\refM, Thm.~1] is wrong e.g.~in the
case when both curves are real conics with non-empty but mutually disjoint
real loci. However, this inaccuracy can be easily corrected and does not
affect the most interesting case when the curves have common real points.
\medskip

\noindent{\bf Definition \defOut.} A component of a $\C$-boundary link $L$
is said to be {\it outer} if it is adjacent to an unbounded component of $B\setminus A$
where $A$ and $B$ are as in the definition of $\C$-boundaries. In particular,
all components of a strong $\C$-boundary are outer.
\medskip

\proclaim{ Proposition \propSum } {\rm(cf.~[\refBR, Prop.~1.4])}
Let $K_1$ and $K_2$ be outer components of $\C$-boundary links
$L_1$ and $L_2$, respectively. Then $L_1\sqcup L_2$ and
$L_1\#L_2=L_1\#_{(K_1,K_2)}L_2$ are $\C$-boundaries.
If, moreover, $L_1$ and $L_2$ are strong $\C$-boundaries, then so are
$L_1\sqcup L_2$ and $L_1\#L_2$, and
$\chi_s(L_1\#L_2)+1=\chi_s(L_1\sqcup L_2)=\chi_s(L_1)+\chi_s(L_2)$.
\endproclaim

\demo{ Proof } For $j=1,2$ let $(A_j,B_j)$ be a realization of $L_j$
as a (strong) $\C$-boundary. By translating $A_1$ and $B_1$ sufficiently far away,
we can achieve that $A_1\cap B_2=A_2\cap B_1=B_1\cap B_2=\varnothing$.
Perturbing $A_1\cup A_2\cup L$ for a suitable line $L$ we can achieve that
$L_j=A\cap\partial B_j$, $j=1,2$, for a smooth projective algebraic curve $A$.
Set $B=B_1\cup B_2\cup T$ where $T$ is a small tubular neighbourhood of an
embedded arc connecting a point in $K_1$ with a point in $K_2$.
Then $(A,B)$ realizes $L_1\sqcup L_2$ (resp. $L_1\#L_2$) if this arc does not
(resp. does) lie on $A$.
Then the required relation on the $\chi_s(\dots)$ easily follows from
Proposition~\propBR.
\qed\enddemo

\proclaim{ Proposition \corBR }
If $L$ is a strong $\C$-boundary and $-L^*$ is a
(not necessarily strong) $\C$-boundary, then $\chi_s(L)=\chi_s^-(L)\ge 1$.
\endproclaim

\demo{ Proof }
% Note that if both $L$ and $L^*$ are strong $\C$-boundaries,
% the result immediately follows from Theorem~\thKM\ combined with Proposition~\propSum.
%
Let $\hat L=L\#_{(K,-K^*)}(-L^*)$
where $-K^*$ is an outer component of $-L^*$, and $K$ is the corresponding component of $L$.
Let $A\cap B$ and $A^*\cap B^*$ be  (strong) $\C$-boundary realizations
of $L$ and $-L^*$, respectively.
The construction in the proof of Proposition~\propSum\ provides a curve
$\hat A$ and a smooth ball $\hat B$ such that $\hat A\cap\hat B$
is $\hat L$ and all the components of $\hat L$ inherited
from $L$ (including $K\#(-K^*)$) lie on the boundary of a single unbounded
component of $\hat A\setminus\hat B$. It is well known (and easy to see)
that such an $\hat L$ bounds a surface
$\Sigma=\Sigma_1\cup\dots\cup\Sigma_r$ in $\hat B$, where $\Sigma_1$ is
a disk bounded by $K\#(-K^*)$ and, for $i\ge2$, \ $\Sigma_i$ is an annulus bounded by
$K_i\sqcup(-K_i^*)$ where $K,K_2,\dots,K_r$ are the
components of $L$. Thus $\chi(\Sigma)=1$.
Let $A'=(\hat A\setminus\hat B)\cup\Sigma$. By construction, $A'$ is connected.
Hence $\chi(A')\le\chi(\hat A)$ by Theorem~\thKM. Again by construction, we have
$\chi(\hat A\cap\hat B) = \chi(A\cap B)+\chi(A^*\cap B^*)-1\le 2\chi_s(L)-1$.
Thus
$$
  0\le\chi(\hat A)-\chi(A')=\chi(\hat A\cap\hat B)-\chi(\Sigma)
   \le (2\chi_s(L)-1)-1.
$$ Finally, $\chi_s(L)=\chi_s^-(L)$ by Proposition \propBR. \qed
\enddemo

\proclaim{ Lemma \lemBBS } Let $L$ be a $\C$-boundary which is
not a strong $\C$-boundary. Then there exists a proper sublink of $L$
which has zero linking number with its complement.
\endproclaim

\demo{ Proof }
Let $L=A\cap B$ as in the definition of a $\C$-boundary.
Without loss of generality we may assume  $A$ to be smooth.
Then $A\setminus B$ has a bounded connected component $A_0$,
because otherwise, for some line $C$,
a perturbation of $A\cup C$ would realize $L$ as a strong $\C$-boundary.
Let $A_1=(A\setminus B)\setminus A_0$, and
let $B'$ be the complement of $B$ in the one-point compactification of $\C^2$.
Then $B'$ is a ball and $A_0$ is disjoint from the closure of $A_1$ in $B'$.
Hence the linking number of $\partial A_0$ and $\partial A_1$
is zero.
\qed\enddemo

%=======================================================

\head \S\sectHopf.
A non-quasipositive $\C$-boundary coming from Wermer's example
\endhead

Let $S^3=\{|z|^2+|w|^2=0\}\subset\C^2$ and let $0<\eps\ll1$. Let $G_f$ be the
graph of the function
$$
       f(z) = \cases 2\eps/\bar z,   & |z|\ge\eps,\\
                     2z/\eps,        & |z|\le\eps.\endcases
$$ which is
endowed with the orientation induced by the projection onto the $z$-axis.
It is easy to check that $L_f:=G_f\cap S^3$ is the link in Figure~\figLink\
where the horizontal circle represents the component of $L_f$ close to the $z$-axis.

\midinsert\centerline{\epsfxsize=40mm\epsfbox{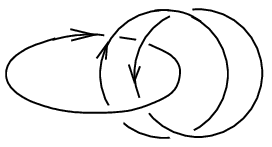}}
\botcaption{Figure~\figHopf} Example of a non-quasipositive $\C$-boundary
\endcaption\endinsert

This link is evidently non-trivial, hence either $L_f$ or its mirror image
$L_{\bar f}$ is not quasipositive by Theorem~\thH.
The mapping $T_f:(z,w)\mapsto(z,w-f(z))$
transforms $G_f$ into the line $\{w=0\}$. Similarly,
$T_{\bar f}(G_{\bar f})=\{w=0\}$. Thus either the pair
$\big(T(B^4),\,\{w=0\}\big)$ or its image under $(z,w)\mapsto(z,\bar w)$ is
a desired example of a ``fancy ball" in $\C^2$ and
an algebraic curve which cuts a non-quasipositive link on its boundary sphere, i.e., either
$L_f$ or $L_{\bar f}$ is a non-quasipositive $\C$-boundary.

So far it is a non-constructive proof of existence, because we do not know yet
which link is not quasipositive (we will establish this later).
However, if we replace $f(z)$ by $f(z+\frac12)+\overline{f(z-\frac12)}$, the
resulting link will be amphicheiral because it is invariant under
$(z,w)\mapsto(-z,-\bar w)$. Thus it is not quasipositive (again by Theorem~\thH)
and it is a $\C$-boundary by the same reasons as above.

Let us show that the link $L_{\bar f}$ is quasipositive (and hence, by Theorem~\thH,
$L_f$ is not). In fact, $L_{\bar f}$ is isotopic to the braid closure
of the $3$-braid $(\sigma_1\sigma_2\sigma_1^{-1})(\sigma_1^{-1}\sigma_2\sigma_1)$.
We do not know whether $L_f$ is a strong $\C$-boundary or not.

The quasipositivity of $L_{\bar f}$ can  also be seen geometrically as follows.
The components of $L_f$ are parametrized, obeying the orientation, as follows
$$
   t\mapsto %\ell_1(t)=
                       (e^{it},\,2\eps e^{it}),\qquad
   t\mapsto %\ell_2(t)=
                       (2\eps e^{-it},\,e^{-it}),\qquad
   t\mapsto %\ell_3(t)=
                       (\tfrac12\eps e^{it},\,e^{it})
$$
(here we have approximated the coefficients up to $O(\eps^2)$).
\if01{

L_f:
   z =         exp( it),  w = 2eps*exp(-it)   ( w = 2*eps / z       )
   z =   2*eps*exp(-it),  w =      exp( it)   ( w = 2*eps / z       )
   z = 1/2*eps*exp( it),  w =      exp(-it)   ( w = 2*conj(z) / eps )
      lk(1,2) = 1    lk(1,3) = -1    lk(2,3) =

L_conj(f):
   z =         exp( it),  w = 2eps*exp( it)   ( w = 2*eps / conj(z) )
   z =   2*eps*exp(-it),  w =      exp(-it)   ( w = 2*eps / conj(z) )
   z = 1/2*eps*exp( it),  w =      exp( it)   ( w = 2*z / eps       )

}\fi
Therefore,
$$
    L_f = S^3\cap
     \Big(\{w=2\eps z\}\cup\{z=2\eps w\}^{\op}\cup\{2z=\eps w\}\Big)
$$
where $\{\dots\}^{\op}$ means the opposite orientation is introduced on a complex line.
%(here we have approximated the coefficients up to $O(\eps^2)$).
Any two triples
of distinct complex lines through the origin are isotopic to each other.
Hence
$$
    L_f \sim S^3\cap
     \Big(\{w=0\}\cup\{w=\eps z\}\cup\{z=\eps w\}^{\op}\Big).    \eqno(\eqA)
    %=S^3\cap\{\bar w(z\bar w-1)=0\}
$$
Thus $L_{\bar f}$ is isotopic to
the image of the right hand side of (\eqA) under $(z,w)\to(z,\bar w)$.
It can be parametrized by
$t\mapsto(e^{it},0)$, $t\mapsto(e^{it}$, $\eps e^{-it})$,
$t\mapsto(\eps e^{-it},\,e^{it})$
but this is (again up to $O(\eps^2)$) a parametrization of
$S^3\cap\{w(zw-\eps)=0\}$. Thus $L_{\bar f}\sim S^3\cap\{w(zw-\eps)=0\}$
is quasipositive.

\medskip\noindent
{\bf Remark \remWermer. }
A third way to see that $L_{\bar f}$ is quasipositve is to observe that
it can be constructed starting from Wermer's example (see  [\refNW, p.~34]),
which consists in exhibiting the function $F(z)=(1+i)\bar z-iz\bar z^2-z^2\bar z^2$
with the following properties: $F'_{\bar z}\ne0$ on the unit disk $\Delta$ and $F|_{\partial\Delta}=0$.
Then the graph of $F$ is totally real, hence it has a small neighbourhood
which is a smooth pseudoconvex ball $B$. One can check that the link cut
by the $z$-axis on $\partial B$ is $L_{\bar f}$; thus $L_{\bar f}$ is quasipositive by [\refO].

%%%%%%%%%%%%%%%%%%%%%%%%%%%%%%%%%%%%%%%%%%%%%%%%%%%%%%%%%%%%%%%%%%
%%%%%%%%%%%%%%%%%%%%%%%%%%%%%%%%%%%%%%%%%%%%%%%%%%%%%%%%%%%%%%%%%%
%%%%%%%%%%%%%%%%%%%%%%%%%%%%%%%%%%%%%%%%%%%%%%%%%%%%%%%%%%%%%%%%%%
%%%%%%%%%%%%%%%%%%%%%%%%%%%%%%%%%%%%%%%%%%%%%%%%%%%%%%%%%%%%%%%%%%

\head   \S\sectFurther. Further examples of $\C$-boundaries cut out by a complex line and their properties
\endhead
It is clear that in \S\sectHopf\ we could take any function $f:\C\to\C$ whose graph $G_f$
is transverse to $S^3$ and cuts a non-trivial link $L_f$ on it. In this
case Hayden's theorem (Theorem \thH) guarantees that either $L_f$ or its mirror image $L_{\bar f}$
is non-quasipositive. This is, however, a very small family of links, which we are
going to describe in this section. All they are iterated torus links, so
an appropriate language to describe them are  {\it EN-diagrams\/}
(called also {\it splice diagrams}), which are certain graphs introduced by Eisenbud and Neumann in [\refEN]. More precisely, considering diagrams obtained one from another by certain simple operations as equivalent (see [\refEN, Thm.~8.1]),
each iterated torus link corresponds to a unique equivalence class of diagrams.

The computation of the iterated torus link structure of $L_f$ from the initial data
is very similar to that in [\refGO] (in both cases, the initial data is an
arrangement of disjoint circles on the plane equipped with some extra information).

So, let $f$, $G_f$, and $L_f$ be as above and let $\pr_1:\C^2\to\C$ be the projection
$(z,w)\mapsto z$. Without loss of generality we assume that $L_f$ is disjoint from the $z$-axis.  %!!!!!!!!
Then $\pr_1(L_f)$ is a disjoint union of smooth embedded
circles $C_1\cup\dots\cup C_n$.
Let $D_1,\dots,D_n$ be the bounded connected components of
$\C\setminus\pr_1(L_f)$ numbered so that $C_j$ is the exterior
component of $\partial D_j$.
We say that $D_j$ is {\it positive} or {\it negative} according to the sign of
$|f(z)|^2+|z|^2-1$ for $z\in D_j$ (thus $\pr_1(G_f\cap B^4)$ is the union of
all negative $D_j$'s).
We endow each $C_j$ with the boundary orientation coming from the adjoining
negative component of $\C\setminus\pr_1(L_f)$ (which is the
orientation induced by the projection of $L_f$).
Let $a_j$ be the
increment of $(\Arg f)/(2\pi)$ along $C_j$.
Then the link $L_f$ is determined by
the following combinatorial data:
the partial order $\prec$ on $\{C_1,\dots,C_n)$ defined by
$C_i\prec C_j$ if $C_i$ lies inside $C_j$, and the numbers $a_j$ corresponding to
the non-maximal $C_j$ with respect to this order.
Such data are realizable if and only if
$\sum_{k\in I(j)} a_k=0$
each time when $D_j$ is positive; here
$I(j)=\{k\mid C_k\subset\partial D_j\}$.

\medskip\noindent{\bf Definition~\defIterTor.}
Let $K$ be a component of an oriented link $L$. Let $p$, $q$, and $d$ be
integers such that gcd$(p,q)=1$ and $d\ge 1$.
We say that $L\cup L'$
(resp. $(L\setminus K)\cup L'$) is the
{\it $(pd,qd)$-cable of $L$ along $K$ with the core retained}
(resp. {\it with the core removed\/}) if, for some tubular neighbourhood $T$
of $K$ disjoint from $L\setminus K$, we have
\roster
\item"$\bullet$" $L'\subset\partial T$, and $L'$ is a disjoint union of knots:
                $L'=K_1\cup\dots\cup K_d$;
\item"$\bullet$" $[K_j]=p[K]$ in $H_1(T)$ and $\lk(K,K_j)=q$ for each $j=1,\dots,d$.
\endroster
An {\it iterated torus link} is a link obtained from the unknot by repeated
cabling of either kind.
\medskip

\noindent{\bf Remark \remIter.}
Reversing the  orientations of some components of an iterated
torus link leads to another iterated torus link. Indeed, reversing
the orientation of a component $K$ is the $(-1,0)$-cable along $K$ with the core
removed.

\proclaim{ Proposition \propIter }
$L_f$ is an iterated torus link.
\endproclaim

\demo{ Proof } This follows from Lemma~\lemIter\ below.
\qed\enddemo

\proclaim{ Lemma \lemIter } Let $L=K_1\cup\dots\cup K_n$ be a link in $S^3=\partial B^4$
such that $\pr_1|_L$ is injective.
Then $L$ is an iterated torus link.
\endproclaim

\demo{ Proof }
We can assume that $L$ is disjoint from the $z$-axis.
Let $K_1,\dots,K_n$ be the components of $L$ and let $C_j=\pr_1(K_j)$.
We shall call the $C_j$   ovals.
Due to Remark~\remIter, on each component we may chose any orientation.
So  we fix on $K_j$ the orientation induced by the counter-clockwise orientation of $C_j$.
Let $a_j$ be the linking number of $K_j$ with the $z$-axis (the
increment of $\Arg F_j/(2\pi)$ where $K_j$ is considered as the graph of a function
$F_j:C_j\to\C$). The link $L$ is determined by $\pr_1(L)$ and the numbers
$a_1,\dots,a_n$ (if $C_j$ is an outermost oval, then $L$ does not depend on this $a_j$
up to isotopy).

We shall prove the lemma for a larger class of links, namely we shall
allow that some components of $L$ are
fibers of $\pr_1$ positively linked with the $z$-axis (in fact
the link does not change if we replace such a component by
a small oval with $a_j=1$).

Without loss of generality we may assume that $\pr_1(L)$ has a single
outermost oval. Otherwise $L$ is a split sum of sublinks each of which
corresponds to an outermost oval and  all the ovals surrounded by it.
If $\pr_1(L)$ consists of a single oval and a point inside it, then
$L$ is the Hopf link which is the $(1,1)$-cable over the unknot.
So it is enough to check that the following operations (i)--(iii)
are cablings
(see the first row in Figure~\figENa).
Let $K$ be a component of $L$ of the form $\pr_1^{-1}(p)$ for a point $p$, and
$D$ be a disk such that $D\cap\pr_1(L)=\{p\}$. The operations are:
\roster
\item"(i)"
    Adding a new component whose projection is $\partial D$ and whose
    linking number with the $z$-axis is any given integer $a$.
\item"(ii)"
    The operation (i) followed by the removal of $K$.
\item"(iii)"
    Replacing $K$ by $\pr_1^{-1}(P)$ where
    $P=\{p_1,\dots,p_k\}\subset D$.
\endroster
Then
 (i) (resp. (ii)) is the $(a,1)$-cabling along $K$
with the core retained (resp. removed), and (iii) is the $(k,0)$-cabling
along $K$ with the core removed.
\qed\enddemo

\midinsert
\centerline{\epsfxsize=120mm\epsfbox{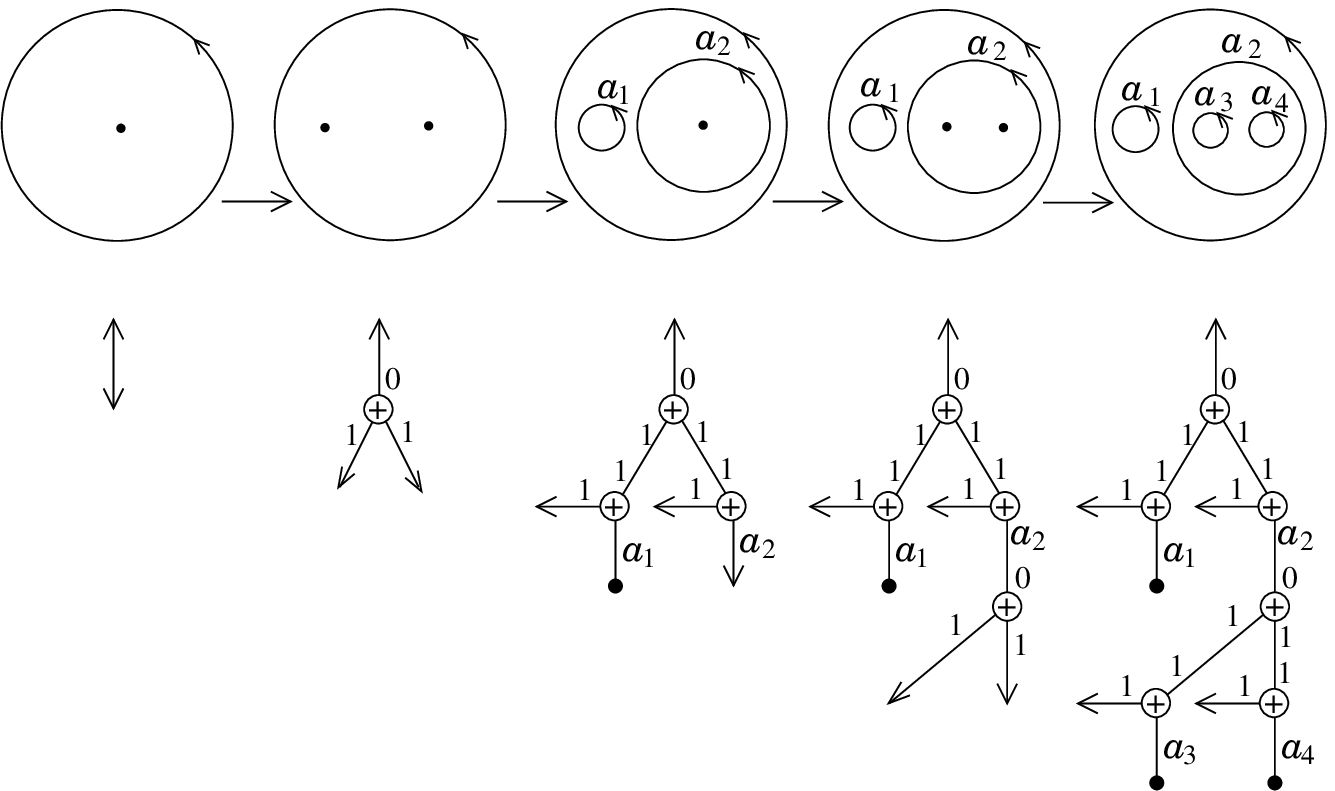}}
\botcaption{ Figure \figENa } Evolution of $\pr_1(L)$ and of the EN-diagram.
\endcaption
\endinsert

The second row in Figure~\figENa\ represents the evolution of the EN-diagram
under the iterated cablings considered in the proof of Lemma~\lemIter.

In Figure~\figENb\ we give two EN-diagrams of a link $L_f$
for the arrangement of ovals and the linking numbers as on the left-hand
side of this figure. The grey area is $\pr_1(G_f\cap B^4)$
(recall that the sum of the linking numbers should be zero along the boundary
of each bounded white (=positive) domain). The left-hand EN-diagram corresponds to the proof of Lemma~\lemIter,
and the right-hand one is obtained from the former  using admissible operations with
EN-diagrams described in [\refEN, Thm.~8.1 (3) and (6)].
In general, such operations can be applied to each piece of an EN-diagram
which corresponds to an annular component of $\pr_1(G_f\setminus B^4)$
(a white annular component for the coloring as in Figure~\figENb).

\midinsert
\centerline{\epsfxsize=110mm\epsfbox{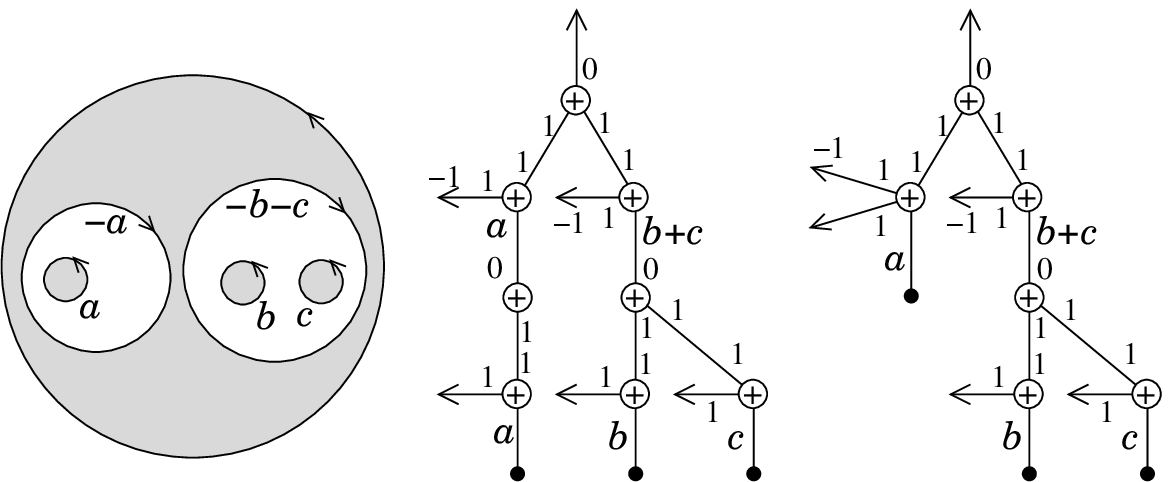}}
\botcaption{ Figure \figENb } $\pr_1(L_f)$ and two EN-diagrams for a $\C$-boundary $L_f$
\endcaption
\endinsert

%=======================================================

\head   \S\sectFiveCr.  Links with at most 5 crossings
\endhead

In this section, for each link admitting a plane projection with at
most five crossings, we determine if it belongs to the
classes $\QP$, $\SCB$ or $\CB$. In Table~1 we give the answers for all such links which
do not have an unknot split component (do not have  the form $L\sqcup O$ where $O$ is
the unknot) but the answers for links of the form
$L\sqcup O\sqcup\dots\sqcup O$ with $\le5$ crossings easily follow.
The implication $L\in\Cal C\Rightarrow L\sqcup O\in\Cal C$
($\Cal C$ is $\QP$, $\SCB$, or $\CB$)
is evident. The reverse implication is true for $\QP$ (see [\refO]) but
we do not know if it takes place in general for the classes $\SCB$ or $\CB$.
However, the nature of our proofs is such that each time when we prove that
$L\not\in\Cal C$ (where $\Cal C$ is the class $\SCB$ or $\CB$), the same arguments can be easily
adapted to prove that $L\sqcup O\sqcup\dots\sqcup O\not\in\Cal C$.

The list of prime links up to 5 crossings is taken from [\refLMknot, \refLMlink]
(but we abbreviate $2_1^2$ to $2_1$). In the second column we give the braid
notation. It serves to identify the link as well as to make evidence of
its quasipositivity (when applicable).
The braid words also help to estimate $\chi_s^-(L)$ from below by using the
observation that if a braid $\beta'$ is obtained from $\beta$ by  replacing some
$\sigma_i^{-1}$ with $\sigma_i$, then  $\chi_s^-(\beta)\ge\chi_s^-(\beta')$.
(In fact, we only need
lower      bounds for $\chi_s^-$ and upper bounds for $\chi_s^-$ for our proofs, and the reader can assume that
Table ~1 contains just these bounds for the Euler characteristics.)
For example,
$\chi_s^-({4_1^2}_-)
=\chi_s^-(\sigma_1\sigma_2^{-1}\sigma_1^{-1}\sigma_1^{-1}\sigma_2^{-1})
\ge\chi_s^-(\sigma_1\sigma_2^{-1}\sigma_1^{-1}\sigma_1\sigma_2^{-1})
=\chi_s^-(\sigma_1\sigma_2^{-2})=2$.
% 1 -2 -1 -1 -2 -> 1 -2 -2 -> -1 -1 => chs=2

We also use Murasugi's inequality in estimates for the links $4_1^2$, ${4_1^2}_-$, and $5_1^2$.

The comments at the end of this section are referred to as ``(a)", ``(b)", etc.  in Table~1.
Almost all proofs are based on general results in \S\S\sectSimpleEx-\sectKM.
Here we present some specific results used in the table.
%We start with proving facts are used in the proofs of the results presented in the table.

\midinsert
\vbox{
\def\TL#1#2#3#4#5#6#7{\noindent
	\hbox to 12mm{\hfill#1\hfill}
	\hbox to 30mm{\hfill#2\hfill}
	\hbox to 18mm{#3\hfill}
	\hbox to 18mm{#4\hfill}
	\hbox to 18mm{#5\hfill}
	\hbox to 10mm{\hfill#6\hfill}
	\hbox to 10mm{\hfill#7\hfill}
}

\hbox to\hsize{\hfill Table 1\;\;\qquad\quad}

\medskip
\noindent\hrule
\smallskip
\TL{$L$}{braid}{$L\in\QP$}{$L\in\SCB$}{$L\in\CB$}{$\chi_s(L)$}{$\chi_s^-(L)$}

\smallskip
\noindent\hrule
\smallskip

\TL{$2_1$}{$\sigma_1^2$}{yes}{yes}{yes}{}{}        % OK

\TL{$2_1^*$}{}{no}{no (a)}{no (b)}{0}{2}           % OK

\medskip

\TL{$3_1$}{$\sigma_1^3$}{yes}{yes}{yes}{}{}        % OK

\TL{$3_1^*$}{}{no}{no (a)}{no (b)}{-1}{1}          % OK

\medskip

\TL{$4_1$}{$(\sigma_1^{-1}\sigma_2)^2$}{no}{no (a)}{no (b)}{-1}{1}  % OK + BR(amphi)

\TL{$4_1^2$}{$\sigma_1^4$}{yes}{yes}{yes}{}{}      % OK

\TL{${4_1^2}^*$}{}{no}{no (a)}{no (b)}{-2}{2}      % OK

\TL{${4_1^2}_-$}{}{no}{no (a)}{no (b)}{0}{2}       % OK
% 1 -2 -1 -1 -2 -> 1 -2 -2 -> -1 -1 => chs=2
% f[3,{-1,2,1,1,2}]  (* {1,0} *) 1 ge 1 + chi => chs le 0
% -1 2 1 1 2 --> -1 1 1 2 = 1 2 => chs ge 0

\TL{${4_1^2}_-^*$}{$\sigma_1^{-1}\sigma_2\sigma_1^2\sigma_2$}{yes}{yes}{yes}{}{}
                                                   % OK

\TL{$2_1\#2_1$}{}{yes}{yes}{yes}{}{}               % OK

\TL{$2_1\#2_1^*$}{}{no (c)}{yes (d)}{yes}{}{}      % OK

\TL{$2_1^*\#2_1^*$}{}{no}{no}{no (f,e)}{-1}{3}   % OK
% f[3,{1,1,2,2}]  (* {2,0} => 1 ge 2 + chi

\TL{$2_1\sqcup 2_1$}{}{yes}{yes}{yes}{}{}          % OK

\TL{$2_1\sqcup 2_1^*$}{}{no}{no (a)}{yes (d)}{0}{2}% OK

\TL{$2_1^*\sqcup 2_1^*$}{}{no}{no}{no (f,e)}{0}{4}     %

\medskip

\TL{$5_1$}{$\sigma_1^5$}{yes}{yes}{yes}{}{}        % OK

\TL{$5_1^*$}{}{no}{no (a)}{no (b)}{-3}{1}          % OK

\TL{$5_2$}{$\sigma_1^2\sigma_2^2\sigma_1\sigma_2^{-1}$}{yes}{yes}{yes}{}{}
                                                   % OK

\TL{$5_2^*$}{}{no}{no (a)}{no (b)}{-1}{1}          % OK
% -1 -1 -2 -2 -1 2 -> 1 -1 2 -2 -1 2 = -1 2 = slice

\TL{$5_1^2$}{$(\sigma_1\sigma_2^{-1})^2\sigma_1$}
                  {no (i)}{yes (j)}{yes}{}{}              % OK

\TL{${5_1^2}^*$}{}{no}{no}{no (f)}{0}{2}            % OK
% <<sm.m; f[3,{1,-2,1,-2,1}] (* {1,0} MT<=> 1 ge 1 + chi

\TL{$3_1\#2_1$}{}{yes}{yes}{yes}{}{}               % OK

\TL{$3_1\#2_1^*$}{}{no (c)}{yes (d)}{yes}{}{}      % OK

\TL{$3_1^*\#2_1$}{}{no}{no}{no (f,g)}{0}{} % OK
% -2 -2 -2 1 1 -> 2 -2 -2 1 1 = -2 1 1 ~ 1 1 => chs=0

\TL{$3_1^*\#2_1^*$}{}{no}{no}{no (f,e)}{-2}{2} % OK
% -1 -1 -1 -2 -2 -> -1 => chs=2

\TL{$3_1\sqcup 2_1$}{}{yes}{yes}{yes}{}{}          % OK

\TL{$3_1\sqcup 2_1^*$}{}{no (c)}{no (a,g)}{yes (d)}{-1}{1}   % OK

\TL{$3_1^*\sqcup 2_1$}{}{no}{no (f,g)}{no (h)}{-1}{1}        %

\TL{$3_1^*\sqcup 2_1^*$}{}{no}{no}{no (f,e)}{-1}{3}    %

\smallskip
\noindent\hrule
\smallskip
} % end vbox
\endinsert

The following theorem is immediate from the
Franks--Williams--Morton inequality [\refFW, \refMo] in combination  with
Proposition~\propBR\ (in [\refLMknot] this result is used in most
cases to prove the non-quasipositivity of knots).

\proclaim{ Theorem~\thFW } {\rm([\refBR, Theorem~3.2]).}
If $L$ is a quasipositive link, then $\ord_v P_L\ge1-\chi_s(L)$,
where $P_L(v,z)$ is the HOMFLY polynomial normalized by
 $P_{\text{unknot}}=1$,
$P_{L_+}=vzP_{L_0}+v^2P_{L_-}$.
\endproclaim

\proclaim{ Proposition \propCaseH }
The link $3_1^*\sqcup 2_1$ is not a $\C$-boundary.
\endproclaim

\demo{ Proof }
Let $L=L_1\sqcup L_2$ where $L_1$ is $3_1^*$ and $L_2$ is $2_1$.
Suppose that $L$ is a $\C$-boundary $A\cap B$.
Let $\Sigma$ be a disjoint union of two surfaces $\Sigma_1\cup\Sigma_2$, and let
$j:(\Sigma,\partial\Sigma)\to(B,\partial B)$ be an immersion with negative
crossings such that $j(\partial\Sigma_i)=L_i$, $i=1,2$.
We may assume that  %??????
 $\Sigma_1$ is a disk and $\Sigma_2$
is an annulus. Let $A'$ be $A\setminus B$ glued with $\Sigma$ along the
boundary, and let us extend $j$ to $A'$, so that
$j(A')=(A\setminus B)\cup j(\Sigma)$.

We have $\chi_s(L)=-1$ and $\chi_s^-(L)=1$
(see Table~1), so that $ \chi_s^-(L)>\chi_s(L)$, and hence      %????????????
$\chi(A)>  \chi(A')$. Thus $A'$ cannot be connected by Theorem \thKM. Hence
 $A'$ % A\setminus B ????????
 has a bounded component $A'_0$ %?????????
  such that
$L_0:=j(A'_0)\cap\partial B$ is a sublink which has zero linking number with
its complement (see Lemma~\lemBBS\ and its proof).
This is possible only %????????
 when $L_0$ is either $L_1$ or $L_2$.
%It is easy to see that $A'_0$ is connected, hence $\chi(A'_0)\le 2$ and,
Since $[j(A'\setminus A'_0)]=[j(A)]$ in $H_2(\CP^2)$, %???????     i chto takoe  j(A)
 by Theorem~\thKM\
we have $\chi(A)\ge\chi(A'\setminus A'_0)$, hence
$$
   \chi(A)+\chi(A'_0)\ge\chi(A')
           = \chi(A\setminus B)+\chi(\Sigma)
           \ge\chi(A)-\chi_s(L)+\chi(\Sigma),
$$
thus $\chi(A'_0) \ge\chi(\Sigma)-\chi_s(L)=2$.
It is easy to see that $A'_0\setminus B$ contains no disk as a component.      Hence $\chi(A'_0)\le 0$, which is
a contradiction.
\qed\enddemo

\proclaim{ Proposition~\propFiveCr } The link $5_1^2$ (see Figure~\figFiveCrA)
is a strong $\C$-boundary.
\endproclaim

\midinsert
\centerline{
   \epsfxsize=35mm\lower-6mm\hbox{\epsfbox{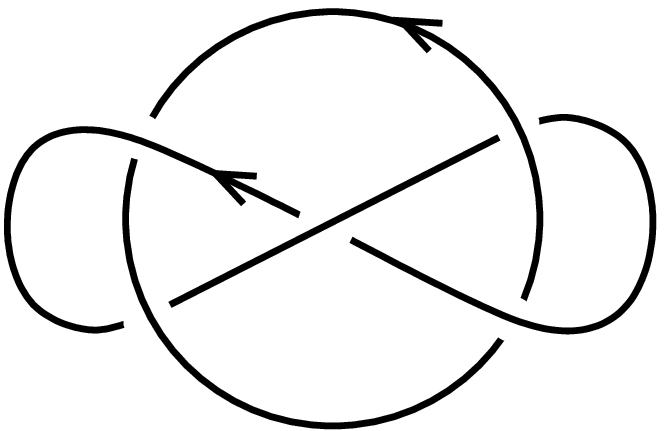}}
   \hskip10mm
   \epsfxsize=42mm\epsfbox{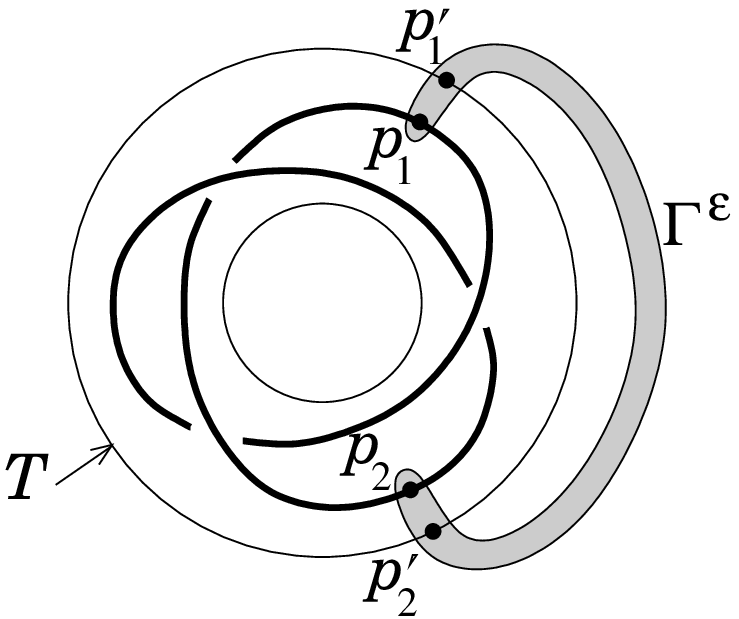}
}
\botcaption{\hskip-15mm Figure~\figFiveCrA. {\rm The link $5_1^2$}
            \hskip 20mm Figure~\figTorus\hskip30mm}
\endcaption
\endinsert

\demo{ Proof }
%Let us fix a small positive $\eps$.
Let $A=\{(z,w)\mid w^2=z^2+z^3\}$. Let $1\ll r\ll R$ and let
$\Delta_r,\Delta_R\subset\C$ be the disks of the respective radii.
Let $U_t=([t,r]\times\Delta_R)\cup\partial(\Delta_r\times\Delta_R)$.
Let $z_1=r\,\exp(\pi i/3)$, $z_2=\bar z_1$, and let $w_j$, $j=1,2$, be
the solution of $w^2=z_j^2+2z_j^3$ such that $\Im w_j>0$, so that
$w_1\approx w_2\approx r^{3/2}i$.
Let $p_j=(z_j,w_j)$ and $p'_j=(z_j,Ri)$.
Let $\gamma$ be an embedded arc in $\Delta_r$ which connects $z_1$ with $z_2$
avoiding the interval $[0,r]$, and let
$\Gamma=[p_1,p'_1]\cup(\gamma\times\{Ri\})\cup[p'_2,p_2]$.
For a set $X\subset\C^2$ and $\eps>0$,
we denote the $\eps$-neighbourhood of $X$ in $\C^2$ by $X^\eps$.
Finally, for $0\ll\delta\ll\eps$, let $B_t$ be a small smoothing of
$(U_t\setminus\Gamma^\eps)^\delta$.

Then $A\cap\partial B_0$ is a strong $\C$-boundary link isotopic to $5_1^2$
in the embedded $3$-sphere $\partial B_0$.
Indeed, we have $U_r=\partial (\Delta_r\times\Delta_R)$, and $A\cap U_r$
is the trefoil knot sitting in the ``vertical'' solid torus
$T=(\partial\Delta_r)\times\Delta_R$; see %?????
 Figure~\figTorus\ where
we represent the piecewise smooth
$3$-sphere $U_r$ via the central projection onto the unit sphere
followed by a suitable stereographic projection onto the $3$-space.

Hence the link $A\cap B_r$ is as on the left-hand side of Figure~\figFiveCr\
(cf.~Theorem~\thB\ and its proof). Consider the family of links
$(B_t,A\cap B_t)$ where $t$ varies from $r$ to $0$. In this deformation,
the link changes only in a small area in its ``inner'' part, namely,
the portion of the link in the sector $-\eta<\Arg z<\eta$ ($0<\eta\ll1$)
of the solid torus $(\partial\Delta_{r-\delta})\times\Delta_R$
is  deformed as $t$ varies.
%by something else during the deformation.
When the parameter $t$ crosses the value $t=\delta$, %????? to zhe samoe ????
 the sphere $\partial B_t$
crosses the double point of $A$ (at the origin), and the link
bifurcates as shown in the middle of Figure~\figFiveCr.
The resulting link is exactly $5_1^2$ (see the right-hand side of Figure~\figFiveCr).
\qed\enddemo

\midinsert
\centerline{\epsfxsize=110mm\epsfbox{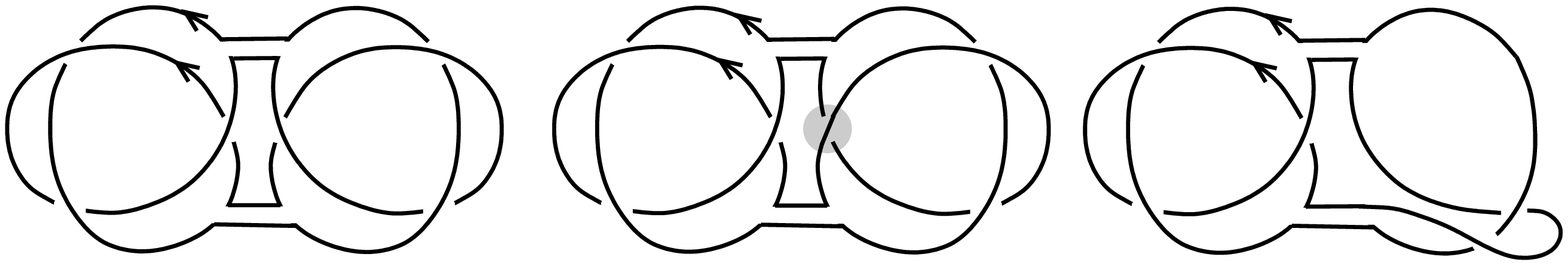}}
\centerline{\hbox to 8mm{}
 $A\cap B_r$ \hskip26.5mm $A\cap B_0$
            \hskip14mm same link redrawn}
\botcaption{Figure~\figFiveCr}
\endcaption
\endinsert

\subhead Comments to Table 1 \endsubhead

(a). Because $\chi_s(L)\ne\chi_s^-(L)$
     (see Proposition~\propBR).

(b). By Lemma~\lemBBS\ combined with the fact that $L\not\in\SCB$.

(c). By Theorem~\thO.

(d). By Theorem~\thA\ applied to the nodal or cuspidal cubic $y^2=ax^2+x^3$
     ($a=0$ or $1$) where $B_0$ is a small ball with center at the origin and
     $B$ a small (for $2_1\#2_1^*$ and $2_1\sqcup 2_1^*$)
     or a large (in the other cases) ball containing $B_0$.
     One easily checks that the resulting link is a strong $\C$-boundary
     in the corresponding cases.

(e). Use that $\chi_s(L)=\chi_s(L^*)$ and
     apply Proposition~\propSum\ to compute $\chi_s(L^*)$.

(f). By Proposition~\corBR.

(g). An embedded surface in a $4$-ball that is bounded by $L$ cannot contain a disk as a component; 
 hence $\chi_s(L)\le0$.
     (In the case of $3_1^*\#2_1$ one can also compute $\chi_s(L)=\chi_s(L^*)$ using Proposition~\propBR.)

(h). See Proposition~\propCaseH.

(i). By Theorem~\thFW.

(j). See Proposition~\propFiveCr.

%=======================================================


\Refs

%\ref\no\refBaader\by S.~Baader
%\paper Slice and Gordian numbers of track knots
%\jour Osaka J. Math. \vol 42 \yr 2005 \pages 257--271 \endref

\ref\no\refBO\by M.~Boileau, S.~Orevkov
\paper Quasipositivit\'e d'une courbe analytique dans une boule pseudo-convexe
\jour C. R. Acad. Sci. Paris, Ser. I \vol 332 \yr 2001 \pages  825--830 \endref

\ref\no\refBR\by M.~Boileau, L.~Rudolph
\paper N\oe uds non concordants \`a un $\C$-bord
\jour Vietnam J. Math. \vol 23 \yr 1995 \pages 13--28 \endref

\ref\no\refEN\by D.~Eisenbud, W.~Neumann
\book Three dimensional Link Theory and Invariants of Plane Curve Singularities
\bookinfo Annals of Math.~Studies 110 \publ Princeton Univ.~Press
\publaddr Princeton NJ \yr 1985 \endref

\ref\no\refFS\by R.~Fintushel, R.~Stern
\paper Immersed spheres in 4-manifolds and the immersed Thom conjecture
\jour Turkish J. of Math. \vol 19 \yr 1995 \pages 145--157 \endref

\ref\no\refFW\by J.~Franks, R.~Williams
\paper Braids and the Jones-Conway polynomial
\jour Trans. Amer. Math. Soc. \vol 303 \yr 1987 \pages 97-108 \endref

% \ref\no\refF\by R.~H.~Fox
% \paper Some problems in  knot theory, Topology of 3-manifolds and related topics
% \inbook Proc. The Univ. of Georgia Institute, 1961 \publ Prentice-Hall
% \publaddr Englewood Cliffs, N. J. \yr 1962 \pages 168-176\endref

\ref\no\refGO\by P.~M.~Gilmer, S.~Yu.~Orevkov
\paper Signatures of real algebraic curves via plumbing diagrams
\jour J.~Knot Theory Ramif. \vol 27 \yr 2018 \issue 3 \pages 1840003, 33 pp \endref

\ref\no\refHay\by K.~Hayden
\paper Minimal braid representatives of quasipositive links
\jour Pac. J. Math. \vol 295 \yr 2018 \pages 421--427 \endref

\ref\no\refKM\by P.~Kronheimer, T.~Mrowka
\paper The genus of embedded surfaces in the projective plane
\jour Math. Res. Letters \vol 1 \yr 1994 \pages 797--808 \endref

\ref\no\refLMknot\by C.~Livingston, A.~H.~Moore
\paper KnotInfo: Table of Knot Invariants \jour \hbox to 35mm{}
http://www.indiana.edu/~knotinfo, June 17, 2020 \endref

\ref\no\refLMlink\by C.~Livingston, A.~H.~Moore
\paper LinkInfo: Table of Link Invariants \jour \hbox to 35mm{}
http://linkinfo.sitehost.iu.edu, June 17, 2020 \endref

\ref\no\refM\by G.~Mikhalkin
\paper Adjunction inequality for real algebraic curves
\jour Math. Res. Letters \vol 4 \yr 1997 \issue 1 \pages 45--52 \endref

\ref\no\refMo\by H.~Morton
\paper Seifert circles and knot polynomials \jour Math. Proc. Camb. Phil.
Soc. \vol 99 \yr 1986 \pages 107--110 \endref

% \ref\no\refMu\by K.~Murasugi
% \paper On the Minkowski unit of slice links
% \jour Trans. Amer. Math. Soc. \vol 114 \yr 1965 \pages 377--383 \endref

\ref\no\refN\by S.~Yu.~Nemirovski
\paper Complex analysis and differential topology on complex surfaces
\jour Russian Math. Surveys \vol 54:4 \yr 1999 \pages 729--752 \endref

\ref\no\refNW\by R.~Nirenberg, R.~O.~Wells
\paper Approximation Theorems on Differentiable Submanifolds of a
Complex Manifold \jour Trans. Amer. Math. Soc.
\vol 142 \yr 1969 \pages 15--35 \endref
% Немировский:
% В самом конце, на с. 34. Но сказано, что пример им указал Вермер.
% Это явная формула для гладкой функции в единичном диске, равной нулю всюду
% на границе и еще в нуле, у которой производная по z с чертой нигде не нуль

\ref\no\refO\by S.~Yu.~Orevkov
\paper Quasipositive links and connected sums
\jour Funk. anal. i prilozh. \vol 54 \yr 2020 \issue 1 \pages 81--86
\lang Russian \transl English transl. \jour Funct. Anal. Appl. \toappear \endref
%\jour arxiv:1906.03454 \endref

\ref\no\refRuTop\by L.~Rudolph
\paper Algebraic functions and closed braids
\jour Topology \vol 22 \yr 1983 \pages 191--201 \endref

\ref\no\refRuFancy\by L.~Rudolph
\paper Plane curves in fancy balls
\jour Enseign. Math. \vol 31 \yr 1985 \pages 81--84 \endref

% \ref\no\refSto\by A.~Stoimenow
% \paper On polynomials and surfaces of variously positive links
% \jour J. Eur. Math. Soc. (JEMS) \vol 7 \yr 2005 \pages 477--509 \endref


\endRefs
\enddocument